# LATTICE POINTS ON THE PLANE AND THE DIOPHANTINE SYSTEM

$$ax + by + cz = d$$
$$a_1 x + b_1 y + c_1 z = d_1$$
$$a_2 x + b_2 y + c_2 z = d_2$$

*By* Konstantine 'Hermes' Zelator

**Lattice Points on the plane** $ax + by + cz = d$ **and**

**the linear diophantine system**

$$a_1x + b_1y + c_1z = d_1$$
$$a_2x + b_2y + c_2z = d_2$$

## Introduction

In almost every introductory number theory book there is a section devoted to the treatment of the linear diophantine equation $ax + by = c$, in two unknowns or variables $x$ and $y$, and with integer coefficients $a, b$ and $c$. The word "diophantine" refers to equations to be solved in the set of integers, usually denoted by Z. Also the same word, "diophantine", has its historic derivation in the name of Diophantos, a Greek mathematician (estimated to have lived in the period 150AD-250AD) who lived in Alexandria, Egypt. Diophantos tackled real life problems whose solution relied on the solution of certain equations to be solved in integers.The books Elementary Number Theory with Applications (see [1]) by Kenneth Rosen and Classical Algebra (see[2]) by William J.Gilbert and Scott A.Vanstone are excellent courses for the student who is unfamiliar with the above linear diophantine equation; they are replete with interesting exercises, not only on this topic but throughout, and they emphasize the computational aspects of number theory in a clear and concise manner. As with every other textbook wiich contains a presentation of the above equation, the complete solution set is presented (parametric solutions in terms of an integer-valued parameter), along with the relevant conditions and facts.

From a geometric perspective when one tries to solve the diophantine equation $ax + by = c$, one has the objective of determining all the lattice points, if any, that lie on the straight line described by the equation $ax + by = c$.In a given coordinate plane with coordinate axes $x$ and $y$, a lattice point $(x_1, y_1)$ is simply a point in which both coordinates $x_1$ and $y_1$ are integers.

Later on in this paper, we will make use of the parametric solution to the above linear diophantine equation. Now let us turn our attention to three-dimensional space equipped with a coordinate system of three mutually perpendicular or orthogonal axes $x, y$ and $z$; and imagine for a moment a lattice points $(x_1, y_1, z_1)$ in space; that would be a point with the real numbers $x_1, y_1, z_1$ actually being integers. The problem of determining the set of all lattice points that lie on a given plane described by the equation $ax + by + cz = d$ (with $a, b, c, d$ being integers) is tantamount to finding all integer solutions to the above equation. This equation often arises in practical problems as well. Only rarely does an introductory number theory book devote a section to this equation. Even more rarely, does it contain a section on the system of two (simultaneous) linear diophantine equations,

$$a_1x + b_1y + c_1z = d_1$$
$$a_2x + b_2y + c_2z = d_2$$

Geometrically speaking we would seek to find the set of all the lattice points which lie on the intersection line of two planes. Of course there may be not a line of intersection (parallel planes) or if there is one, it may not contain any lattice points. This is then the aim of this article: to present (in a step-by-step manner) the parametric solutions to the above three-variable diophantine equation as well as the system of the two three-variable diophantine equations; along with interesting computational examples.

The organization of this paper in sequential sections is as follows:

| |
|---|
| 1. Terminology, 2. Review: Solving the diophantine equation $ax + by = c$, 3. The Four Formulas,4. Examples, 5. Explanations and Proofs, 6. Chart of Special Cases |

## 1. Terminology

A notation of the form $(k, \ell)$ can be interpreted three ways: as a point in the coordinate plane (with axes $x$ and $y$); secondly as the greatest common divisor of two integers $k$ and $\ell$.And thirdly, as a solution to the diophantine equation $ax + by = c$ (here, we exempt the open interval notation).

Similarly, a notation of the form $(k, \ell, r)$ can have three meanings: as a point in space (equipped with an $x, y, z$ coordinate system), as the greatest common divisor of three integers $k, \ell, r$ ,or (third meaning) as a solution to the linear diophantine equation $ax + by + cz = d$ or the diophantine system

$$\left.\begin{array}{l} a_1x + b_1y + c_1z = d_1 \\ a_2x + b_2y + c_2z = d_2 \end{array}\right\}$$

It will always be clear from the text whether for example $(k, \ell, r)$ refers to the greatest common divisor of the integers $k, \ell, r$; or to a point in space or solution to the above tree variable diophantine equation or system. Note that the difference between the meanings of the triple $(k, \ell, r)$ as a point or as solution to the above equation (or system) is only semantic; for all practical purposes we may regard the frist and the third meanings as indistinguishable. Thus the only distinction that needs to be made clear is that between $(k,\ell,r)$ as a point (or solution) and as greatest common divisor of $k,\ell,r$, we will be writing $(k,\ell,r)=\delta$ or some other letter.
With regard to the notation $x \in$ Z which may appear on occasion throughout this paper, it simply means "$x$ belongs to the set Z" or "$x$ is a member of the set Z (in plain language, "$x$ is an integer").

Finally, the conjunctive "$\wedge$" may be used between statements, which simply means "and".

And one more clarification. If we state that "the integer k divides the integer $\ell$", we will mean that k exactly divides $\ell$; that is, $k=t\cdot\ell$, for some $t\in$ Z. The same meaning we will attribute to the expression " k is a divisor of $\ell$". The notation "|" can also be employed "$k \mid \ell$" means "k divides $\ell$".

**2. Review: Solving the diophantine equation** $ax+by=c$

First assume that both a and b are nonzero integers.The special case (and the geometric interpretation of) $ab=0$ (i.e. $(a=0, b\neq 0)$ , $(a\neq 0, b=0)$ or $(a=b=0)$ is discussed in the special cases section.(Section 6)

The graph described by the equation $ax+by=c$ is slant (straight) line whose slope is equal to $-\dfrac{a}{b}$. The equation $ax+by=c$ has solutions in the set of integers Z if, and only if the greatest common divisor $(a,b)$ is also a divisor of $c : \delta=(a,b)\mid c$.

If $\delta$ does not divide c then the set of solutions is the empty set (no solutions) We have, $ax+by=c \Leftrightarrow \left(\dfrac{a}{\delta}\right)x+\left(\dfrac{b}{\delta}\right)y=\dfrac{c}{d}$; the last diophantine equation is equivalent to the original one; that is, they have the same solution set. The two reduced integers $\dfrac{a}{\delta}$ and $\dfrac{b}{\delta}$ are relatively prime or coprime, another way of saying that their greatest common divisor is equal to 1: $\left(\dfrac{a}{\delta},\dfrac{b}{\delta}\right)=1$. If $(x_1,y_1)$ is a particular solution to the above linear diophantine equation then all the solutions are given by

Solution set S: $(x,y)=\left(x_1-\left(\dfrac{b}{\delta}\right)m, y_1+\left(\dfrac{a}{\delta}\right)m\right)$

(For further details, proof etc, see either of the two references, [1] or [2])
Where m is an integer parameter; in other words m can take any integer value; we can simply write $m\in$Z.Thus, for each specific integer value of m a new particular solution is produced.

So for $m=0,(x,y)=(x_1,y_1)$. Thus, if we know one specific or particular solution, all the rest can be found. But how does one determine a particular solution $(x_1,y_1)$? Sometimes it is easy to find one by inspection or by a short trial and error process, especially if the coefficients a and b are small in absolutely value. But there is fail proof procedure that works in all cases, although it may several steps in order to arrive at a particular solution $(x_1,y_1)$.Essentially this process is what is known as the Euclidean Algorithm. We can explain how it works without becoming unduly theoretical. For simplicity let us put $a_1=\frac{a}{\delta},b_1=\frac{b}{\delta},c_1=\frac{c}{\delta};(a_1,b_1)=1$.We have the equation $a_1x+b_1y=c_1$. If $|a_1|=1$ or $|b_1|=1$, we are done; for example say $|a_1|=1$ (i.e. $a_1=1$ or -1) and suppose $a_1=-1$.We solve for x to obtain $x=-c_1+b_1y$. So the solution set S is given by S: $(x,y)=(-c_1+b_1\cdot m,m), m\in$ Z. A similar picture emerges in the cases $a_1=1, b_1=1$ or $b_1=-1$. Now suppose that $|a_1|>1$ and $|b_1|>1$. Because of the condition $(a_1,b_1)=1$ it follows that $|a_1|$ and $|b_1|$ must be distinct positive integers: either $|a_1|>|b_1|$ or $|a_1|<|b_1|$. For the sake of discussion suppose that $|a_1|<|b_1|$: we solve the equation for x (if it were $|a_1|>|b_1|$, we would solve for y):

$$a_1x+b_1y=c_1 \Leftrightarrow x=\frac{c_1-b_1y}{a_1}$$

We perform the divisions $c_1:a_1$ and $b_1:a_1$; we have

$$b_1=a_1q_1+r_1, c_1=a_1Q_1+R_1,$$

Where $q_1,Q_1$ are the quotients of the divisions and $r_1,R_1$ the two remainders.

As we know form the division theorem, $0\le r_1<|a_1|$ and $0\le R_1<|a_1|$. We find that,

$$x=\frac{(a_1q_1+r_1)-(a_1Q_1+R_1)y}{a_1}$$

$$x=q_1-Q_1y+\frac{r_1-R_1y}{a_1}$$

Obviously for x to be an integer, the ratio $\frac{r_1-R_1y}{a_1}$ must be an integer; say

$$\frac{r_1-R_1y}{a_1}=z_1 \Leftrightarrow r_1=a_1z_1+R_1y.$$

Note that in the last equation, $R_1 = |R_1| < |a_1|$. Thus by setting $y = z_0$,we can repeat the process by solving for the variable $z_0$ (which is one that has the smaller, in absolute value, coefficient). The Euclidean Algorithm ensures us that the process must terminate after a finite number of steps by obtaining, in the last stage, an equation of the form $kz_{n+1} + z_n = \mu$, which has the integer solution $(\mu,0)$.By back substitution we can trace a particular solution $(x_1, y_1)$ to equation $a_1x + b_1y = c_1$. Let us illustrate this with an example.

Consider the linear diophantine equation,

102x + 140y = 318

We have a = 102, b = 140 , c = 318; $\delta$ = (a,b) = (102 , 140) = 2;

$a_1 = \frac{a}{\delta} = 51$, $b_1 = \frac{140}{2} = 170$, and $c_1 = \frac{c}{\delta} = 159$. The original equation is equivalent to 51x + 70y = 159.

Obviously $|51| = 51 < |70| = 70$, so we solve for x:

$$x = \frac{159 - 70y}{51}$$

The divisions 159:51 and 70:51 yield 159=(3)(51)+6 and 70=(1)(51)+19.
Therefore,

$$x = \frac{[(3)(51)+6]-[(1)(51)+19]y}{51};$$

$$x = 3 - y + \frac{6-19y}{51}$$

Next we set $y = z_0$ and $z_1 = \frac{6-19z_0}{51}$, a linear diophantine equation with unknowns $z_0$ and $z_1$; we obtain $51z_1 + 19z_0 = 6; |19| = 19 < |51| = 51$,so we solve for $z_0$:

$$z_0 = \frac{6-51z_1}{19} = \frac{6-(2.19+13)z_1}{19} \Leftrightarrow z_0 = -2z_1 + \frac{6-13z_1}{19}$$

Again, we set $z_2 = \frac{6-13z_1}{19}$, a linear diophantine equation in the unknowns $z_1$ and $z_2$: we have $19z_2 + 13z_1 = 6$; note that $|13| = 13 < |19| = 16$ and so solve for $z_1$:

$$z_1 = \frac{6-19z_2}{13} = \frac{6-(1.13+6)z_2}{13}$$

$$z_1 = -z_2 + \frac{6(1-z_2)}{13}$$

We could continue for a few more steps till we get an equation of the form $kz_{n+1} + z_n = \mu$ (which has the integer solution $(\mu,0)$ ); but it is not really necessary, for the last diophantine equation has an obvious particular solution: $z_2 = 1$ and $z_1 = -1$. Tracing back we find $y = z_0 = -2(-1)+1 = 3$ and $x = 3-3-1 = -1$

Thus $(x_1, y_1) = (-1,3)$ is a particular solution to the diophantine equation 51x + 70y = 159 and hence the original one. All the solutions are given by

$x = x_1 - b_1 m$ , $y = y_1 + a_1 m$ ;

**Solution** set S(*x*,*y*)=(-1-70*m*,3+51*m*), $m \in \mathbb{Z}$

In the Examples Section, in all the examples (except for one) the particular solutions are easily found by inspection. However, we wanted to give an illustration to the reader of how this process works in less than obvious cases.

### 3. The Four Formulas

In formulas 1,2 and 3 below assume that the coefficients a,b,c are nonzero (for the cases in which abc=0, see the Special Cases Section).

**Formula 1**

| |
|---|
| If at least on of the positive integers $\lvert a\rvert, \lvert b\rvert, \lvert c\rvert$ is equal to 1 (i.e. $a = \pm 1$ *or* $b = \pm 1$ *or* $c = \pm 1$), all the solutions to the diophantine equation $ax + by + cz = d$ can be obtained by merely solving for the corresponding variable. For example if a=1, solution set $S : (x, y, z) = (d - bm - cn, m, n)$, where m and n are integer- value parameters.<br>While in the case b=-1,<br>Solution set $S : (x, y, z) = (m, am + cn - d, n)$, *with* $m, n \in \mathbb{Z}$ |

The next formula pertains to the case in which $\lvert a\rvert > 1, \lvert b\rvert > 1, \lvert c\rvert > 1$ and at least one of the greatest common divisors (a,b),(b,c),(a,c) is equal to 1. First observe that the diophantine equation $ax + by + cz = d$ will have solutions if, and only if, the number $\delta = (a,b,c)$ is a divisor of d. Otherwise, if $\delta$ does not divide d, the equation has no integer solutions. If $\delta \mid d$ , the given diophantine equation is equivalent to $\left(\frac{a}{\delta}\right)x+\left(\frac{b}{\delta}\right)y+\left(\frac{c}{\delta}\right)z=\frac{d}{\delta}$; and $\left(\frac{a}{\delta},\frac{b}{\delta},\frac{c}{\delta}\right)=1$

Because of the above reduction, in Formulas 2 and 3, we start with an equation which is already reduced;

**(*a*,*b*,*c*)=1**

**Formula 2**

Assume (*a, b, c*) = $1, |a|, |b|, |c| > 1$ and at least one of (*a,b*), (b,c), (a,c) being equal to 1. If (*a,b*) = 1, all the solutions to the diophantine equation $ax+by+cz=d$ are given by $x=(d-cm)x_1-bn,\ y=(d-cm)y_1+an, z=m$, with $m,n$ being integer-valued parameters, and ( $x_1, y_1$ ) being a solution to the linear diophantine equation $ax+by=1$

Likewise, if $(b,c)=1$, all the solution can be expressed by the two-parameter formulas $x=m, y=(d-am)y_1-cn, z=(d-am)z_1+bn, with\,(x_1, y_1)$ being a solution to the equation $by+cz=1$.

Similarly, if $(a,c)=1$, all the solutions can be given by $x=(d-bm)x_1-cn$, $y=m, z=(d-bm)z_1+an$ with $(x_1, y_1)$ being a solution to $ax+cz=1$

The next formula, Formula 3, deals with the cases in which not only each of the three coefficient is greater than 1 in absolute value; but also each of the greatest common divisors $(a,b),(b,c),(a,c)$ is greater than 1 (so we can not use Formula 2).

There is a fact from number theory that is used in Formula 3, namely $((a,b),c)=(a,b,c)$ for any nonzero integers $a,b,c$. We leave this as an exercise for the reader to show. In our case,

**(*a*,*b*,*c*)=1= ((*a*,*b*),*c*)**

**Formula 3**

Assume $|a|,|b|,|c|>1$ and $(a,b),(b,c),(a,c)>1$, and let $(a,b)=\delta$

All the solutions to the diophantine equation

$ax+by+cz=d$ are given by the formulas,

$$\left.\begin{aligned} x &= (t_1-cm)x_1-\left(\frac{b}{\delta}\right)n \\ y &= (t_1-cm)y_1+\left(\frac{a}{\delta}\right)n \\ z &= z_1+\delta m \end{aligned}\right\}$$

Where $m,n$ are integer-valued parameters, $(x_1,y_1)$ is a solution to the linear diophantine equation $\left(\frac{a}{\delta}\right)x+\left(\frac{b}{\delta}\right)y=1$; and $(t_1,z_1)$ is a solution to the linear diophantine equation $\delta t+cz=d$.

To finish this section, we stat Formula 4 which pertains to the diophantine system mentioned in the introduction.

**Formula 4**

Assume that at most one of the six integers $a_1,b_1,c_1,a_2,b_2,c_2$ is zero and consider the linear diophantine system

$$\left.\begin{aligned} a_1x+b_1y+c_1z &= d_1 \\ a_2x+b_2y+c_2z &= d_2 \end{aligned}\right\}$$

Where $d_1$ and $d_2$ are also integers.

i. If $a_1=\varepsilon a_2 \wedge b_1=\varepsilon b_2 \wedge c_1=\varepsilon c_2$ where $\varepsilon=1$ or -1, then the above equation has solutions if, and only if, $d_1=\varepsilon d_2$; that is, if $d_1\neq\varepsilon d_2$, the solution set is the empty set; while for $d_1=\varepsilon d_2$, the solution set to the above system is equal to the solution set of the linear diophantine equation $a_1x+b_1y+c_1z=d_1$ which can be solved with the aid of Formulas 1,2 or 3; provided that $(a_1,b_1,c_1)$ is a divisor of $d_1$ (otherwise there are no solutions); if one of $a_1,b_1$ or $c_1$ happens to be zero, refer to the Special Cases section.

ii. Assume that at least one of the conditions $a_1 \neq \varepsilon a_2$, $b_1 \neq \varepsilon b_2$, $c_1 \neq \varepsilon c_2$ holds true as well as $(a_1, b_1, c_1) = (a_2, b_2, c_2) = 1$.

Let $D_1, D_2, D_3, D$ be the determinants respectively, of the 2x2 matrices

$$\begin{bmatrix} a_1 & b_1 \\ a_2 & b_2 \end{bmatrix}, \begin{bmatrix} b_1 & c_1 \\ b_2 & c_2 \end{bmatrix}, \begin{bmatrix} a_1 & c_1 \\ a_2 & c_2 \end{bmatrix}, \begin{bmatrix} d_1 & c_1 \\ d_2 & c_2 \end{bmatrix}:$$

$D_1 = a_1 b_2 - a_2 b_1$, $D_2 = b_1 c_2 - b_2 c_1$, $D_3 = a_1 c_2 - a_2 c_1$ and $D = d_1 c_2 - d_2 c_1$.

$D_2$ and $D_3$ can not both be zero (see Explanations and Proofs section) and so we can define,

$$D_{23} = (D_2, D_3) and\, \delta = \left( \frac{c_1 . D_1}{D_{23}}, c_1 \right)$$

Under the condition $c_1 \neq 0$ (see notes 1 and 2 below).

The above diophantine system will have solutions (i.e. a nonempty solution set) only if both conditions, stated below, are satisfied:

**Condition 1:** The integer $D_{23}$ divides (or is a divisor of) the integer D

**Condition 2:** $\delta$ is a divisor of the integer $d_1 - a_1 x_1 - b_1 y_1$, for any solution ($x_1, y_1$) to the linear diophantine equation $D_3 x + D_2 y = D$ (see note 3).

If both conditions are met, the solution set to the above system can be described by,

$$x = x_1 - \left( \frac{D_2}{D_{23}} \right) m_1 + \left( \frac{c_1}{\delta} \right) \left( \frac{D_2}{D_{23}} \right) \lambda,$$

$$y = y_1 + \left( \frac{D_3}{D_{23}} \right) m_1 - \left( \frac{c_1}{\delta} \right) \lambda, z = z_1 + \frac{1}{\delta} \left( \frac{c_1 d_1}{D_{23}} \right) \lambda,$$

where $\lambda$ is and integer-valued parameter and $(m_1, z_1)$ is a solution to the linear diophantine equation,

$$\left( \frac{c_1 d_1}{D_{23}} \right) m + c_1 z = d_1 - a x_1 - b y_1$$

(which will have solutions, by condition 1).

**Note 1:** Note that the number-fraction $\frac{c_1 d_1}{D_{23}}$ is indeed an integer: an easy computation shows that,

$$-a_1 \left( \frac{D_2}{D_{23}} \right) + b_1 \left( \frac{D_3}{D_{23}} \right) = \frac{c_1 D_1}{D_{23}}$$

Obviously , both $\frac{D_2}{D_{23}}$ and $\frac{D_3}{D_{23}}$ are integers since $D_{23}$ is the greatest common divisor of $D_2$ and $D_3$.Similarly, the numbers $\frac{c_1}{\delta}$ and $\left( \frac{c_1 D_1}{D_{23}} \right) \frac{1}{\delta}$ are both integers by the very definition of the integer $\delta$.

**Note 2:** We may assume $c_1 \neq 0$, since if it were $c_1 = 0$; by hypothesis, we know that at most one of the coefficient is zero and so by necessity we would have $c_2 \neq 0$.We can simply rename the coefficients: $a_1, b_1, c_1$ become $a_2, b_2, c_2$ and conversely.

**Note 3:** With regard to Condition 1: It will either hold for every solution $(x_1, y_1)$or for none; this will become evident in the Explanations and Proof section.

When one is faced with a specific linear diophantine equation (or system) of two, three, or more variables, one may solve it without applying any systematic general method or formulas; but if one follows a different path for solving the very same system one may end up with a different set of parametric formulas describing the same solution set. Assuming that no error occurred in the process, both sets of parametric formulas would be correct; the would be equivalent: one could, in a finite number of algebraic steps derive one from the other. We make this remark in light of the possibility that some readers may want to experiment and improvise with specific linear diophantine equations or systems and discover this fact in their own experience.

## 4. Examples

Below we present seven examples.

- **Example 1**

  i. Find the solutions to the diophantine equation $x - 3y - 4z = 0$
  ii. Determine the lattice points that lie on the plane described by the above equation and which also lie in the interior of or on the boundary of (sphere) the ball whose center is the origin (0,0,0) and whose radius is R=2.

**Solution**

i) We are in the simplest of cases, namely those dealt with by Fromula 1. We solve for x to obtain x = 3y + 4z ; solution set $S:(x,y,z)=(3m+4n,m,n)$, where $m,n$ are integer-valued parameters.

ii) For a point $(x,y,z)$to lie inside or on the given ball, it is necessary and sufficient that its distance from the center not exceed R=2; since the center is the point (0,0,0), we must have,

$$x^2 + y^2 + z^2 \leq 4$$

It follows that

$$x^2, y^2, z^2 \leq 4 \Leftrightarrow (|x| \leq 2 \wedge |y| \leq 2 \wedge |z| \leq 2) \Leftrightarrow$$
$$\Leftrightarrow (-2 \leq x \leq 2, -2 \leq y \leq 2), -2 \leq z \leq 2$$

These conditions are necessary but not sufficient for a point $(x, y, z)$ to lie in or on the given ball.The locus that the last three inequalities describe, is the set of all points in space which lie inside or on the rectangular box whose boundary surface is a cube centered at (0,0,0) and whose edges have length 4.
4. The given ball is inscribed in that box. Applying part (i), we arrive at the necessary conditions, $(-2 \le 3m + 4n \le 2 \wedge -2 \le m \le 2 \wedge -2 \le n \le 2)$.

A quick search can be done systematically by picking a value of m(m=-2,-1,0,1,2) and checking to see whether there are values of n satisfying the first and the third inequalities. Seven points are produces: (-2,-2,1) , (2,-2,2), (1,-1,1) , (0,0,0) , (-1,1,-1) , (-2,2,-2) and (2,2,-1).

However only three lie in or the given ball (actually all three lie inside the ball):

(1,-1,1) , (0,0,0) , (-1,1,-1)

- **Example 2**

In how many ways can a person pay the amount of 80 cents using only dimes, nickels or quarters?

**Solution**

If x is the number of dimes, y the number of nickels, and z the number of quarters, we must have 10x + 5y + 25z = 80 $\Leftrightarrow$ 2x + y + 5z = 16.

We are seeking the number of nonnegative solutions to the last diophantine equation; we must have $0 \le x, y, z$ which lead to the constraints $0 \le x \le 8 \wedge 0 \le 16 \wedge 0 \le z \le 3$. This is a Formula 1 case, solving the equation yields x = m , y=16 – 2m – 5n , z = n, where m and n are integer-valued parameters. Applying the constrains implies,

$$0 \le m \le 8 \wedge 0 \le 2m + 5n \le 16 \wedge 0 \le n \le 3$$

A search produces exactly twenty solutions; the triples (0,16,0),(0,11,1),(2,6,2),(0,1,3),(1,14,0),(1,9,1),(1,4,2),(2,12,0),(2,7,1),(2,2,2),(3,10,0) ,(3,5,1),(3,0,2),(4,8,0),(4,3,1),(5,6,0),(5,1,1,),(6,4,0),(7,2,0), and (8,0,0)

Therefore there are exactly twenty ways of paying the given amount.

- **Example 3**

i) Find the integer solutions to the equation 2x + 3y + 7z = 23
ii) Find those lattice points that lie on the plane described by the equation of part (i) and in the interior space bounded by or on the cube which is described by the inequalities $-3 \le x \le 3 \wedge -3 \le y \le 3 \wedge -3 \le z \le 3$.

**Solution**

i. Since (2,3) = (3,7) = (2,7) = 1, we can apply Formula 2 with more than one choice; we apply it with (*a*,*b*) = (2,3) = 1; in our case $a$ = 2 , b = 3 , $c$ = 7, and $d$ = 23. All the solutions can be described by $x$ = (23-7m) $x_1$ - 3n, $y = (23-7m)y_1 + 2n, z = m, where\ m, n \in$ Z are the two parameters and $(x_1, y_1)$ is an integer solution to the equation 2$x$ + 3$y$ =1.

By inspection, $(x_1, y_1) = (-1,1)$ is solution. Hence,
Solution set $S : (x, y, z) = (-23 + 7m - 3n, 23 - 7m + 2n, m)$.

ii. We apply the constrains to the solutions we found:
$(-3 \le -23 + 7m - 3n \le 3 \wedge -3 \le 23 - 7m + 2n \le 3 \wedge -3 \le m \le 3) \Leftrightarrow$
$\Leftrightarrow (20 \le 7m - 3n \le 26 \wedge 20 \le 7m - 2n \le 26 \wedge -3 \le m \le 3)$

A search shows that only $m$=2 , m=1 and $m$ = 3 yield values of $n$ that satisfy the first and second inequalities. Specifically for $m$=2 we obtain $n$ = -3 and $n$ = -4; while for $m$ = 1 the values of $n$ are only $n$ = -6; and for $m$=3 those values of $n$ are $n$= -1 and $n$ = 0. Going back to the parametric formulas for $x$,$y$, and $z$ we obtain exactly five points.

($x$,$y$,$z$) = (0,3,2) , (3,1,2) , (2,4,1) , ( 1,0,3) , (-2,2,3)

The third point in the list lies in the interior space bounded by the cube, whereas the other four lie on the cube itself.

- **Example 4**

Find the solutions to the diophantine equation 6$x$ – 15$y$ + 10z = 4.

**Solution**

We have $a$ = 6 , $b$ = -15 , $c$ = 10 , $d$ = 4 thus ($a$,$b$,$c$) = 1. Since $|a|, |b|, |c| > 1$, $(a, b) = 3 > 1, (b, c) = 5 > 1,$ and $(a, c) = 2 > 1$, the relevant formula to use is Formula 3.

First we need to find an integer solution $(x_1, y_1)$ to the linear equation
$\left(\frac{a}{\delta}\right)x + \left(\frac{b}{\delta}\right) = 1$, where $\delta = (a, b) = 3$.

By inspection the equation 2$x$ – 5$y$ = 1 has a solution $(x_1, y_1) = (3,1)$. Next we find and integer solution $(t_1, z_1) = (-2,1)$ is a solution. Applying the parametric formulas in Formula 3 we obtain,

Solution set S:($x$,$y$,$z$) = (-6 - 30$m$ + 5$n$, -2 – 10$m$ + 2$n$, 1 + 3$m$), with $m$ ,$n$ $\in$ Z

Check: 6 (-6 – 30$m$ + 5$n$) -15 (-2 – 10$m$ + 2$n$) + 10(1 + 3$m$) = 4

- **Example 5**

Find the solutions to the linear diophantine system $\left.\begin{array}{l} 6x - 4y + 3z = 30 \\ 3x + 6y - 2z = 25 \end{array}\right\}$

**Solution**

We have $a_1 = 6, b_1 = -4, c_1 = 3, d_1 = 30, a_2 = 3, b_2 = 6, c_2 = -2, d_2 = 25$.
Note that $(a_1, b_1, c_1) = 1 = (a_1, b_2, c_2)$ and all the coefficients $a_1, b_1, c_1, a_2, b_2, c_2$ are nonzero. Also, obviously, the requirement that $a_1 \neq \pm b_1$ or $a_2 \neq \pm b_2$ or $a_3 \neq \pm b_3$ is satisfied, so we can apply Formula 4, part (ii). We compute the four determinants to find $D_1 = 48, D_2 = -10, D_3 = -21, D = -135$. Also,

$D_{23} = (D_2, D_3) = (-10, -21) = 1$ and $\delta = \left(\frac{c_1 \cdot D_1}{D_{23}}, c_1\right)$ = ((3) (48),3)=3.

We check the two conditions:

**Condition 1**: $D_{23} = 1$ obviously divides D = -135. So this condition is also satisfied.

**Condition 2**: First we find a solution $(x_1, y_1)$ to the equation

-21x – 10y = - 135 $\Leftrightarrow$ 21x + 10y = 135; $(x_1, y_1)$=(5,3) is a solution.

We compute, $d_1 - a_1x_1 - b_1y_1 = 12$ which is divisible by $\delta$ = 3; so the condition is satisfied.

Thus, the system has solutions. To find them, we must find a solution $(m_1, z_1)$ to the linear equation,

$$\left(\frac{c_1 \cdot D_1}{D_{23}}\right)m + c_1 z = d_1 - (a_1x_1 + b_1y_1);$$

144m + 3z = 12 $\Leftrightarrow$ 48m + z = 4; an obvious solution is $(m_1, z_1) = (0,4)$. Applying this together with the other information we already have to the formulas in Formula 4. we find,

Solution set $S: (x, y, z) = (5 - 10\lambda, 3 + 21\lambda, 4 + 48\lambda); \lambda \in$ Z
Check: $6(5 - 10\lambda) - 4(3 + 21\lambda) + 3(4 + 48\lambda) = 30$
And $3(5 - 10\lambda) + 6\,(3 + 21\lambda) - 2\,(4 + 48\lambda) = 25$

- **Example 6**

Solve the diophantine system $\left.\begin{aligned} 13x + 11z &= 123 \\ -5y + 7z &= 4 \end{aligned}\right\}$

**Solution**

Observe that this system falls in the category of the Special Cases section since two of the six coefficients $a_1, b_1, c_1, a_2, b_2, c_2$ are zero, namely $b_1$ and $a_2$. We can solve this system by first finding the general solution to the second equation and substituting (for y) into the first equation. By inspection $y_1 = z_1 = 2$ is a solution to the second equation. The general solution can be given by,

$$y = 2 - 7\lambda, z = 2 - 5\lambda; \lambda \in \mathbb{Z}$$

Substituting into the first equation gives

$$13x + 11(2 - 5\lambda) = 123 \Leftrightarrow 13x - 55\lambda = 101$$

To find a solution $(x_1, \lambda_1)$ to the last equation, we solve for x: we have

$$x = \frac{55\lambda + 101}{13} = 4\lambda + 7 + \frac{3\lambda + 10}{13}$$

Let t = $\frac{3\lambda + 10}{13} \Leftrightarrow \lambda = \frac{13t - 10}{3} = 4t - 3 + \frac{t-1}{3}$

An obvious solution to the last equation is $t_1 = 1 = \lambda_1$

Thus $x_1 = 4\lambda_1 + 7 + \frac{3\lambda_1 + 10}{13} = 4 + 7 + \frac{13}{13} = 12$

We conclude that $(x_1, \lambda_1) = (12,1)$ is a solution to the linear diophantine equation $13x - 55\lambda = 101$. Hence all integer solutions to the last equation are given by x=12-(-55) m and $\lambda$=1 + 13m. We substitute for $\lambda$ in the previous equations for y and z to find.

Solution set S: (x,y,z) = (12+55m , 5-91m , -3-65m).
Check: 13(12+55m) + 11 (-3-65m) = 123
and -5 (-5-91m) +7(-3-65m) = 4

- **Example 7**

Let p be a positive integer.

i) Describe the set of all triples (x,y,z) with the two properties,

1. The sum of the three integers x,y,z is equal to p, and
2. The integers 7x , 5y, 3z (in that order) are successive terms of an arithmetic progression.

ii) What condition must p satisfy in order that x,y,z be positive integers?
iii) Find the smallest value of p for which x,y,z are positive integers.
iv) What conditions must p satisfy in order that x,y,z be the side lengths of a triangle?
v) For p = 85 find those triples (x,y,z) with x,y,z>0. Which of those triples correspond to a triangle with x,y,z being the side lengths?

**Solution**

i) if the numbers 7x , 5y , 3z are consecutive terms of an arithmetic progression, the middle term must be the average of the other two. Thus we obtain the diophantine system

$$\left\{\begin{array}{l} x+y+z=p \\ \\ 2(5y)=7x+3z \end{array}\right\} \Leftrightarrow \left\{\begin{array}{l} z=p-(x+y) \\ \\ 7x-10y+3z=0 \end{array}\right\} \Leftrightarrow$$

$$\Leftrightarrow \left\{\begin{array}{l} z=p-(x+y) \\ 7x-10y+3[p-(x+y)]=0 \\ \end{array}\right\} \Leftrightarrow \left\{\begin{array}{l} z=p-(x+y) \\ -4x+13y=3p \\ \end{array}\right\}$$

In the last system, we solve the second equation for x to obtain $x=3y+\dfrac{y-3p}{4}$ ; clearly, $(x_1,y_1)=(9p,3p)$ is a solution, and so all the solutions to the second linear diophantine equation are given by $x=9p-13m$ and $y=3p-4m$, where m is an integer-valued parameter. Substituting back in the first equation (of the last system above) for z we find that,

Solution set $S:(x,y,z)=(9p-13m,3p-4m,-11p+17m), m\in \mathbb{Z}$

ii) Setting $(x>0\wedge y>0\wedge z>0)\Leftrightarrow$
$\Leftrightarrow (9p-13m>0\wedge 3p-4m>0\wedge -11p+17m>0)\Leftrightarrow$
$\Leftrightarrow \left(m<\dfrac{9p}{13}\wedge m<\dfrac{3p}{4}\wedge m<\dfrac{11p}{4}\right)\Leftrightarrow$
$\Leftrightarrow \left(\text{since } p>0 \text{ and } \dfrac{11}{17}<\dfrac{9}{13}<\dfrac{3}{4}\right)\dfrac{11p}{17}<m<\dfrac{9p}{13}$

is the necessary and sufficient condition for *x,y* and z to be positive integers. What this really says is that in order for *x,y* and z to be positive integers it is necessary and sufficient that the open interval $l_p=\left(\dfrac{11p}{17},\dfrac{9p}{13}\right)$(here we use standard precalculus/ calculus notation) contain at least one integer *m*. For each such value of the parameter *m*, a positive integer triple (*x,y,z*) will be generated.

iii) An easy calculation show that $p = 3$ is the smallest value of the (positive) integer p for which the interval $l_p$ contains an integer. For $p = 3$ , $l_p = l_3 = \left(\frac{33}{17}, \frac{27}{13}\right); m = 2$ is the only integer that falls in $l_3$. For p = 3 and m = 2 the formulas part (i) produce (*x*,*y*,*z*) = (1,1,1).

iv) If (*x*,*y*,*z*) corresponds to a triangle, the three triangle inequalities must be satisfied: $(x + y > z \wedge y + z > x \wedge x + z > y)$. Note that these three inequalities alone imply $x > 0 \wedge y > 0 \wedge z > 0$( for example, add the first two inequalities member wise to obtain $x > 0$ ).

This means that we should obtain an interval (that depends on *p*) which must contain at least one integer, and which is a subinterval of the interval obtained in part (ii). Indeed, using the above three triangles inequalities and parametric formulas for *x*,*y*,*z* ( in part (i) we now have,

$$(12p - 17m < -11p - 17m \wedge -8p + 13m > 9p - 13m \wedge -2p + 4m > 3p - 4m) \Leftrightarrow$$

$$\Leftrightarrow \left( m > \frac{23p}{34} \wedge m > \frac{17p}{26} \wedge m > \frac{5p}{8} \right) \Leftrightarrow$$

$$\Leftrightarrow \left( \sin ce\ p > 0\ and\ \frac{5}{8} < \frac{17}{20} < \frac{23}{34} \right) \frac{17p}{26} < m \frac{23p}{34} .$$

We see that the necessary and sufficient condition for (*x*,*y*,*z*) to correspond to a triangle is that the open interval $J_p = \left(\frac{17p}{26}, \frac{23p}{34}\right)$ contain and integer. Note as expected, that since $\frac{11}{17} < \frac{17}{26} < \frac{23}{34} < \frac{9}{13}$, the interval $J_p$ lies entirely within $l_p$(the interval in part (ii)).

v) For $p = 85 \Leftrightarrow l_{85} = \left(\frac{935}{17}, \frac{765}{13}\right)$ and $J_{85} = \left(\frac{1445}{26}, \frac{1955}{34}\right)$. Also since $\frac{935}{17} = 55, \frac{765}{13} \approx 58.846, \frac{1445}{26} \approx 55.569, \frac{1955}{34} \approx 57.5$, we see that the internal $l_{85}$ contains three integers; the integers m = 56,57,58; whereas the interval $J_{85}$ contains two integers, namely m = 56,57.

Substituting for *p* = 85 and *m* = 56,57,58 in the formulas in part (i) we find that,

(*x*,*y*,*z*) = (37,31,17) , (24,27,34) , (11,23,51)

The first two triples correspond to a triangle, but the third one does not since 11 + 23 = 34 is in fact less than 51; this is as expected by virtue of the fact that the integer 58 does not lie within the interval $J_{85}$. Also note that in part (iii) where we found that for $p = 3$, there is one positive integer triple $(x,y,z)$, namely $(x,y,z) = (1,1,1)$; we see that $p = 3$ is also the smallest value of $p$ for which a triple $(x,y,z)$ corresponds to triangle; which for $p = 3$ is the equilateral triangle of side lengthe 1.

5. **Explanations and Proofs**

1. **Of Formula 1**

Not much to show here, one simply solves for the variable whose coefficient is 1 or -1.

2. **Of formula 2**

We only need explain the derivation of Formula 2 in the case $(a,b) = 1$. Since $(a,b) = 1$, the linear diophantine equation $ax + by = 1$ has solutions, let $(x_1, y_1)$ be a particular solution. We have $ax_1 + by_1 = 1$. We put z = m , where m can be any integer; so that,

$$ax_1 + by_1 = 1 \Leftrightarrow a[(d-cm)x_1] + b[(d-cm)y_1] = d - cm,$$

which shows that the pair $((d-cm)x_1, (d-cm)y_1)$ is a solution to the linear diophantine equation ax + by =d – cm. Consequently, all the integer solutions of the last equation are given by $x = (d-cm)x_1 - bn, y = (d-cm)y_1 + an, n \in$ Z. But the initial diophantine equation $ax + by + cz = d$ is obviously equivalent to $(ax + by = d - cz \wedge z = m, m \in Z)$. Hence all the solutions to the initial equation are given by $x = (d-cm)x_1 - bn, y = (d-cm)y_1 + an, z = m$, where m,n are integer-valued parameters. We are done.

3. **Of Formula 3**

The diophantine equation $ax + by + cz = d$ is equivalent to $ax + by = d - cz$. Since $\delta$ is the greatest common divisor of a and $b$, we must have $a = \delta\left(\frac{a}{\delta}\right)$ and $b = \delta\left(\frac{b}{\delta}\right)$, where the integers $\frac{a}{\delta}$ and $\frac{b}{\delta}$ are relatively prime; that is, $\left(\frac{a}{\delta}, \frac{b}{\delta}\right) = 1$. We see that,

$$ax + by = d - cz \Leftrightarrow \delta\left(\frac{a}{\delta}\right)x + \delta\left(\frac{b}{\delta}\right) = d - cz \Leftrightarrow$$

$$\Leftrightarrow \delta\left[\left(\frac{a}{\delta}\right)x + \left(\frac{b}{\delta}\right)y\right] = d - cz$$

Clearly the last equation will have solutions if, and only if, $\delta$ is divisor of $d$-$cz$. Evidently the original equation is equivalent to the system of diophantine equations,

$$\left.\begin{array}{l}\left(\frac{a}{\delta}\right)x+\left(\frac{b}{\delta}\right)y=t\\ \delta t=d-cz\end{array}\right\}\Leftrightarrow\left.\begin{array}{l}\left(\frac{a}{\delta}\right)x+\left(\frac{b}{\delta}\right)y=t\\ \delta t+cz=d\end{array}\right\}\quad\begin{array}{l}(1)\\(2)\end{array}$$

the variables being *x,y,z* and *t*.

Observe that each equation in the last system has solutions: the first has solutions since $\left(\frac{a}{\delta},\frac{b}{\delta}\right)=1$. The second equation has solutions as well by virtue of $(\delta,c)=((a,b),c)=(a,b,c)=1$, which is true by hypothesis. If $(x_1,y_1)$ is an integer solution $\left(\frac{a}{\delta}\right)x+\left(\frac{b}{\delta}\right)y=1$, then obviously $(tx_1,ty_1),t\in$ Z, is a solution to equation (1). Furthermore, let $(t_1,z_1)$ be a particular solution to (2); then all the integer solutions to equation (2) are given by

$$t=t_1-cm, z=z_1+\delta m, m\in \text{Z}$$

But for a given integer value of t, $(tx_1,ty_1)$ is particular solution to equation (1); therefore all the integer solutions to equation (1), while (2) also holds true, can be described by $x=tx_1-\left(\frac{b}{\delta}\right)n, y=ty_1+\left(\frac{a}{\delta}\right)n$; where $n\in$ Z and $t=t_1-cm$. We conclude that the solution set to the original equation can be described by,

$$x=(t_1-cm)x_1-\left(\frac{b}{\delta}\right)n, y=(t-cm)y_1+\left(\frac{a}{\delta}\right)n, z=z_1+\delta m;$$

Where m,n are integer-valued parameters, $(x_1,y_1)$ is a particular solution to the linear diophantine equation $\left(\frac{a}{\delta}\right)x+\left(\frac{b}{\delta}\right)y=1$ and $(t_1,z_1)$ a particular solution to $\delta t+cz=d$.

4.. **Of Formula 4**

i) This part is obvious, it is self-explanatory.

ii) First note that the integers $D_2$ and $D_3$ can not be both zero: for if that were the case we would have $a_1c_2-a_2c_1=0=b_1c_2-b_2c_1$; then none of $a_1,a_2,b_1,b_2,c_1,c_2$ could be zero, since if one them were zero, another would also have to be zero (for example, if $a_1=0$, then $a_1c_2-a_2c_1=0\Leftrightarrow a_2c_1=0$, which is not allowed by hypothesis, since at most one of the six coefficients can be zero. Thus, all six coefficients would be nonzero and so we would have $\frac{a_1}{a_2} = \frac{b_1}{b_2} = \frac{c_1}{c_2}$; but also, $\frac{c_1}{c_2} = \frac{c_1'(c_1, c_2)}{c_2'(c_{1,}c_2)}$, where $c_1'$ and $c_2'$ are relatively prime integers (i.e. $(c_1', c_2') = 1$) and $(c_1, c_2)$ is the greatest common divisor of $c_1$ and $c_2$; thus, $\frac{a_1}{a_2} = \frac{b_1}{b_2} = \frac{c_1'}{c_2'}$, which implies (see remark 2)
$(a_1 = \mu c_1', a_2 = \mu c_2', b_1 = \lambda c_1', b_2 = \lambda c_2')$,
for some integers $\mu$ and $\lambda$; we see that $a_1 = \mu c_1', b_1 = \lambda c_1', c = c_1'(c_1, c_2)$ and in view of the assumption that greatest common divisor of $a_1, b_1$ and $c_1$ is equal to 1, we conclude $c_1' = 1$ or -1; Likewise the hypothesis $(a_2, b_2, c_2) = 1$ implies $c_2' = 1$ or -1. Obviously there are four combinations of values of $c_1'$ and $c_2'$ but in all cases $\frac{c_1'}{c_2'} = \varepsilon$, where $\varepsilon = 1$ or -1. Consequently,

$$\frac{a_1}{a_2} = \frac{b_1}{b_2} = \frac{c_1}{c_2} = \frac{c_1'}{c_2'} = \varepsilon \Rightarrow (a_1 = \varepsilon a_2 \wedge b_1 = \varepsilon b_2 \wedge c_1 = \varepsilon c_2),$$

Contrary to the assumption of part (ii).

It is now clear that at least one of $D_2, D_3$ must be nonzero and therefore their greatest common divisor $D_{23}$ exists and is, of course, a positive integer. We have,

$$\left\{ \begin{matrix} a_1x + b_1y + c_1z = d_1 \\ a_2x + b_2y + c_2z = d_2 \end{matrix} \right\} \underset{\text{since } c_1 \neq 0}{\Leftrightarrow} \left\{ \begin{matrix} a_1x + b_1y + c_1z = d_1 \\ -c_1a_2x - c_1b_2y - c_1c_2z = -c_1d_2 \end{matrix} \right\} \Leftrightarrow$$

$$\Leftrightarrow (c_2 \text{ times } 1^{\text{st}} \text{ eq.plus } 2^{\text{nd}}) \quad \left\{ \begin{matrix} a_1x + b_1y + c_1z = d_1 & (3) \\ D_3x + D_2y = D & (4) \end{matrix} \right\}$$

We see that in order for equation (4) to have integer solutions it is necessary and sufficient that $D_{23} = (D_2, D_3)$ be a divisor of D, in other words, condition 1 must hold true.

If $(x_1, y_1)$ is an integer solution to equation (4), then all integer solutions to (4) are given by $x = x_1 - \left(\frac{D_2}{D_3}\right)m, y = y_1 + \left(\frac{D_3}{D_{23}}\right)m, m \in Z^+$ (5)

Substituting for $x$ and $y$ in equation (4) yields, $\left(-a_1\left(\frac{D_2}{D_{23}}\right) + b_1\left(\frac{D_3}{D_{23}}\right)m\right) + c_1 z = d_1 - (a_1 x_1 + b_1 y_1)$, and since the quantity in the brackets in equal to $\frac{c_1 D_1}{D_{23}}$ (see note 1) we arrive at

$$\left(\frac{c_1 D_1}{D_{23}}\right)m + c_1 z = d_1 - (a_1 x_1 + b_1 y_1) \quad (6)$$

It is clear that the original diophantine system is equivalent to the pair of equations (5) and (6) in the variables *m* and *z*. For equation (6) to have integer solutions, it is necessary and sufficient that $\delta = \left(\frac{c_1 D_1}{D_{23}}, c_1\right)$ be a divisor of the integer $d_1 - (a_1 x_1 + b_1 y)$; in other words, condition 2 must hold true. The integer solutions to (6) given by, $m = z_1 - \left(\frac{c_1}{\delta}\right)\lambda,\ z = m_1 + \frac{1}{\delta}\left(\frac{c_1 D_1}{D_{23}}\right)\lambda,$ where $\lambda$ is an integer-valued parameter and $(m_1, z_1)$ a particular solution to (6). Substituting for *m* and *z* in (5) we find

$$x = x_1 - \left(\frac{D_2}{D_{23}}\right)z_1 + \left(\frac{D_2}{D_{23}}\right)\lambda,$$

$$y = y_1 + \left(\frac{D_3}{D_{23}}\right)z_1 + \left(\frac{D_3}{D_{23}}\right)\lambda,$$

$$z = m_1 + \frac{1}{\delta}\left(\frac{c_1 D_1}{D_{23}}\right)\lambda, \qquad \square$$

**Remark 1**. Note the assumption $(a_1, b_1, c_1) = 1 = (a_2, b_2, c_2)$. Obviously if we are given a system of two linear diophantine equations in three variables, the first thing to do would be check whether the greatest common devisor of the coefficients of the unknowns in each equation divides the corresponding constant $(d_1 \; or \; d_2)$ on the other side of the equation. If either equation fails the test, system has no integer solutions. If on the other hand, both equation pass the test, we reduce each equation the (by dividing both sides with $d_i$) so that we have $(a_i, b_i, c_i) = 1$ for $i = 1,2$

**Remark 2**. In the explanations of formula 4 we made use of the following fact from number theory: if four non zero integers $\alpha, \beta, \gamma, \delta$ satisfy the conditions $\frac{\alpha}{\beta} = \frac{\gamma}{\delta}$ and $(\gamma, \delta) = 1$ then $(\alpha = k\gamma \wedge \beta = k\,\delta)$ for some nonzero integer *k*. Indeed $\frac{\alpha}{\beta} = \frac{\gamma}{\delta}$ implies $\alpha\delta = \beta\gamma$; since $\gamma$ is relatively prime to $\delta$ and $\gamma$ divides the product $\alpha\delta$, it must divide $\alpha$;this fact from number theory is typically covered in the first 2 or 3 weeks of an introductory course in elementary number theory. We have $\alpha = k\gamma$ for some non zero integer $k$; and from $\alpha\delta = \beta\gamma \Rightarrow (k\gamma)\delta = \beta\gamma \Rightarrow$ (since $\gamma \neq 0$) $\beta = k\delta$.

**6. Chart of Special Cases**

A. *The linear diophantine equation* $ax + by = c$ with $ab = 0$.

i) If $a = 0$ *and* $b \neq 0$, then the equation has no solution if *b* is not a divisor of *c*. If ,on the other hand *b* is a devisor of *c* then the solution set consists of all pairs of the form $\left(m, \frac{b}{c}\right)$, where $m \in$ Z(all lattice points on the horizontal line $y = \frac{c}{b}$).

ii) If $a \neq 0$ *and* $b = 0$, there are no solutions if *a* is not a divisor of *c*; otherwise (if *a* divides *c*), solution set $S : (x, y) = (\frac{c}{a}, m), m \in$ Z (all lattice points on the vertical line $x = \frac{c}{a}$)

iii) If $a = b = 0$ there are no solutions unless $c = 0$, in which case the solution set consist of all lattice points on the plane. Namely, $S : (x, y) = (m, n), m, n \in$ Z

B. *The linear diophantine equation* $ax + by + cz = d$ with $abc = 0$.

Let $[a \quad b \quad c]$ be the matrix of the coefficients. In each group below, we only present the solution for the first (as viewed from left to right) matrix (of the coefficients) in the group. The others are treated similarly.

i) *Group* 1: Exactly one of the coefficients is zero. Possible matrices are $\begin{bmatrix}0 & b & c\end{bmatrix}, \begin{bmatrix}a & 0 & c\end{bmatrix}, \begin{bmatrix}a & b & 0\end{bmatrix}$.

$\begin{bmatrix}0 & b & c\end{bmatrix}$: The equation describes a plane which is parallel to or contains the *x*-axis. If $(b,c)$ does not divide *d*, there are no solutions. If, on the other hand $\delta = (b,c)$ divides *d*, the solution set is given by

$S:(x,y,z) = \left(m, y_1 - \left(\frac{c}{\delta}\right)n, z_1 + \left(\frac{b}{\delta}\right)n\right); m,n \in Z$ and $(y_1, z_1)$ a solution to $\left(\frac{b}{\delta}\right)y + \left(\frac{c}{\delta}\right)z = \frac{d}{\delta}$.

ii) *Group* 2: Two of the coefficients are zero. The third non zero
$\begin{bmatrix}a & 0 & 0\end{bmatrix}, \begin{bmatrix}0 & b & 0\end{bmatrix}, \begin{bmatrix}0 & 0 & c\end{bmatrix}$

$\begin{bmatrix}a & 0 & 0\end{bmatrix}$: If *a* is not a divisor of *d* there are no solutions. If *a* does divide *d*, solution set $S:(x,y,z) = \left(\frac{d}{a}, m, n\right); m,n \in Z$. The equation describes a plane parallel or coincident to the *yz*- plane.

iii) *Group* 3: $\begin{bmatrix}0 & 0 & 0\end{bmatrix}$

*There are no solutions unless* $d = 0$, *in which case* $S:(x,y,z) = (m,n,\lambda), m,n,\lambda \in Z$ *(all lattice points in space).*

C. The linear diophantine system $\left.\begin{matrix} a_1x + b_1y + c_1z = d_1 \\ a_2x + b_2y + c_2z = d_2 \end{matrix}\right\}$, with at least two of the coefficients $a_1, b_1, c_1, a_2, b_2, c_2$ being zero. Matrix of coefficients is $\begin{bmatrix} a_1 & b_1 & c_1 \\ a_2 & b_2 & c_2 \end{bmatrix}$

**1. Exactly two of the six coefficients are zero**

i) Group1: Two zeros in a column.

$\begin{bmatrix} 0 & b_1 & c_1 \\ 0 & b_2 & c_2 \end{bmatrix}, \begin{bmatrix} a_1 & 0 & c_1 \\ a_2 & 0 & c_2 \end{bmatrix}, \begin{bmatrix} a_1 & b_1 & 0 \\ a_2 & b_2 & 0 \end{bmatrix}$

$\begin{bmatrix} 0 & b_1 & c_1 \\ 0 & b_2 & c_2 \end{bmatrix}$: Each equation describes a plane that is parallel to or contains the *x*-axis. If $b_1c_2 - b_2c_1 = 0$ or equivalently $\frac{b_1}{b_2} = \frac{c_1}{c_2}$; then the system has no solutions unless $\frac{b_1}{b_2} = \frac{c_1}{c_2} = \frac{d_1}{d_2}$, in which case the problemis reduced to the single equation $ax + by + cz = d$ with $a_1 = 0$ (see part B, case (i)).

If $b_1c_2 - b_2c_1 \neq 0$, there are no solutions unless the integer $(b_1c_2 - b_2c_1)$ is a common divisor of the integers $(d_1c_2 - d_2c_1)$ and $(d_2b_1 - b_2d_1)$;if that is the case, solution set $S:(x,y,z) = \left( m, \frac{d_1c_2 - d_2c_1}{b_1c_2 - b_2c_1}, \frac{d_2b_1 - b_2c_1}{b_1c_2 - b_2c_1} \right), m \in \mathbb{Z}$.

ii) *Group* 2: Two zeros in a row.

$$\begin{bmatrix} 0 & 0 & c_1 \\ a_2 & b_2 & c_2 \end{bmatrix} \begin{bmatrix} a_1 & 0 & 0 \\ a_2 & b_2 & c_2 \end{bmatrix} \begin{bmatrix} 0 & b_1 & 0 \\ a_2 & b_2 & c_2 \end{bmatrix} \begin{bmatrix} a_1 & b_1 & c_1 \\ 0 & 0 & c_2 \end{bmatrix} \begin{bmatrix} a_1 & b_1 & c_1 \\ a_2 & 0 & 0 \end{bmatrix} \begin{bmatrix} a_1 & b_1 & c_1 \\ 0 & b_2 & 0 \end{bmatrix}$$

$\begin{bmatrix} 0 & 0 & c_1 \\ a_2 & b_2 & c_2 \end{bmatrix}$: The first equation describes a plane parallel to of coincident with *xy* plane. If $c_1$ is a not a divisor of $d_1$, the system has no integer solutions. If $c_1$ is a divisor of $d_1$, the system will have solutions only if $\delta = (a_2, b_2)$ is a divisor of the integer $d_2 - c_2\left(\frac{d_1}{c_1}\right)$; in that case all the solutions are given by, $S:(x,y,z) = \left( x_1 - \left(\frac{b_2}{\delta}\right)m, y_1 + \left(\frac{a_2}{\delta}\right)m, \frac{d_1}{c_1} \right)$, where $(x_1, y_1)$ is a solution to $\left(\frac{a_2}{\delta}\right)x + \left(\frac{b_2}{\delta}\right)y = \frac{d_2 - c_2\left(\frac{d_1}{c_1}\right)}{\delta}$.

iii) *Group* 3: Two zeros not both on the same row or column.

$$\begin{bmatrix} 0 & b_1 & c_1 \\ a_2 & 0 & c_2 \end{bmatrix} \begin{bmatrix} 0 & b_1 & c_1 \\ a_2 & b_2 & 0 \end{bmatrix} \begin{bmatrix} a_1 & 0 & c_1 \\ 0 & b_2 & c_2 \end{bmatrix} \begin{bmatrix} a_1 & 0 & c_1 \\ a_2 & b_2 & 0 \end{bmatrix} \begin{bmatrix} a_1 & b_1 & 0 \\ 0 & b_2 & c_2 \end{bmatrix} \begin{bmatrix} a_1 & b_1 & 0 \\ a_2 & 0 & c_2 \end{bmatrix}$$

$\begin{bmatrix} 0 & b_1 & c_1 \\ a_2 & 0 & c_2 \end{bmatrix}$: If $\delta_1 = (b_1, c_1)$ is not a divisor of $d_1$, the system will have no solutions .If $\delta_1$ is a divisor of $d_1$ the system will have solutions only if $\delta_2 = \left( a_2, c_2\left(\frac{b_1}{\delta_1}\right) \right)$ is a divisor of the integer $d_2 - c_2z_1$, where $(y_1, z_1)$ can be any particular solution to $b_1y + c_1z = d_1$; if that is the case, the solution set to the system can be described by,
$S: x = x_1 - \frac{1}{\delta_2}\left( c_2\frac{b_1}{\delta_1} \right), y = y_1 - m_1\left(\frac{c_1}{\delta_1}\right) - \left(\frac{c_1}{\delta_1}\right)\left(\frac{a_2}{\delta_2}\right)\lambda$ ,and

$z = z_1 + m_1\left(\frac{b_1}{\delta_1}\right) - \left(\frac{b_1}{\delta_1}\right)\left(\frac{a_2}{\delta_2}\right)\lambda$ ,where $(x_1, m_1)$ is a solution to

$a_2 x + \left(c_2\left(\frac{b_1}{\delta_1}\right)\right)m = d_2 - c_2 z_1$ and $\lambda$ is an integer-valued parameter.

**2. Exactly three of the six coefficients are zero**

i) The two of the three zeros are on the same column and two among them lie on the same row. Possible matrices:

$$\begin{bmatrix} 0 & 0 & c_1 \\ 0 & b_2 & c_2 \end{bmatrix}\begin{bmatrix} 0 & b_1 & 0 \\ 0 & b_2 & c_2 \end{bmatrix}\begin{bmatrix} 0 & b_1 & c_1 \\ 0 & b_2 & 0 \end{bmatrix}\begin{bmatrix} 0 & b_1 & c_1 \\ 0 & 0 & c_2 \end{bmatrix}\begin{bmatrix} 0 & 0 & c_1 \\ a_2 & 0 & c_2 \end{bmatrix}\begin{bmatrix} a_1 & 0 & 0 \\ a_2 & 0 & c_2 \end{bmatrix}$$

$$\begin{bmatrix} a_1 & 0 & c_1 \\ a_2 & 0 & 0 \end{bmatrix}\begin{bmatrix} a_1 & 0 & c_1 \\ 0 & 0 & c_2 \end{bmatrix}\begin{bmatrix} 0 & b_1 & 0 \\ a_2 & b_2 & 0 \end{bmatrix}\begin{bmatrix} a_1 & 0 & 0 \\ a_2 & b_2 & 0 \end{bmatrix}\begin{bmatrix} a_1 & b_1 & 0 \\ a_2 & 0 & 0 \end{bmatrix}\begin{bmatrix} a_1 & b_1 & 0 \\ 0 & b_2 & 0 \end{bmatrix}$$

$\begin{bmatrix} 0 & 0 & c_1 \\ 0 & b_2 & c_2 \end{bmatrix}$: The first equation describes a plane parallel to or coincident with *xy*-plane; the second equation describes a plane parallel to or containing the *x*- axis. The system will have no solutions only if $c_1$ not a divisor of $d_1$. If on the other hand, $c_1 \mid d_1$, the system will have solutions only if $b_2$ is a divisor of the integer $d_2 - \left(\frac{d_1}{c_1}\right)c_2$. In that case,

Solution set $S:(x, y, z) = \left( m, \frac{d_2 - \frac{d_1}{c_1}c_2}{b_2}, \frac{d_1}{c_1} \right)$

ii) Two among the three zeros lies on the same row, but no two lie on the same column. Possible matrices:

$$\begin{bmatrix} 0 & 0 & c_1 \\ a_2 & b_2 & 0 \end{bmatrix}\begin{bmatrix} 0 & b_1 & 0 \\ a_2 & 0 & c_2 \end{bmatrix}\begin{bmatrix} a_1 & 0 & 0 \\ 0 & b_2 & c_2 \end{bmatrix}\begin{bmatrix} a_1 & b_1 & 0 \\ 0 & 0 & c_2 \end{bmatrix}\begin{bmatrix} a_1 & 0 & c_1 \\ 0 & b_2 & 0 \end{bmatrix}\begin{bmatrix} 0 & b_1 & c_1 \\ a_2 & 0 & 0 \end{bmatrix}$$

$\begin{bmatrix} 0 & 0 & c_1 \\ a_2 & b_2 & 0 \end{bmatrix}$: The first equation describes a plane parallel to or coincident with the *xy*-plane; the second equation describes a plane parallel to or containing the *z*-axis. If $c_1$ is not a divisor of $d_1$, there are no solutions. If on the other hand, $c_1 \mid d_1$, the system will have solutions only if $\delta = (a_2, b_2) \mid d_2$. If that is the case, the solution set can be described by $S:(x,y,z)=\left(x_1-\left(\frac{b_2}{\delta}\right)m,\ y_1+\left(\frac{a_2}{\delta}\right)m,\frac{d_1}{c_1}\right)$, $m\in\mathbb{Z}$; where $(x_1,y_1)$ is a solution to $\left(\frac{a_2}{\delta}\right)x+\left(\frac{b_2}{\delta}\right)y=\frac{d_2}{\delta}$.

iii) Three zeros on the row.

$\begin{bmatrix}0 & 0 & 0\\ a_2 & b_2 & c_2\end{bmatrix}\begin{bmatrix}a_1 & b_1 & c_1\\ 0 & 0 & 0\end{bmatrix}$

$\begin{bmatrix}0 & 0 & 0\\ a_2 & b_2 & c_2\end{bmatrix}$: If $d_1\neq 0$ there are no solutions, while if $d_1=0$, the system is reduced to $a_2x+b_2y+c_2z=d_2$ with $a_2b_2c_2\neq 0$, which has already been covered by Formulas 1,2 or 3.

3. **Exactly four of the six coefficients are zero**

i) *Group* 1: $\begin{bmatrix}0 & 0 & c_1\\ 0 & 0 & c_2\end{bmatrix}\begin{bmatrix}0 & b_1 & 0\\ 0 & b_2 & 0\end{bmatrix}\begin{bmatrix}a_1 & 0 & 0\\ a_2 & 0 & 0\end{bmatrix}$

$\begin{bmatrix}0 & 0 & c_1\\ 0 & 0 & c_2\end{bmatrix}$: There are solutions only if $c_1\mid d_1$ and $c_2\mid d_2$; if so, there are solutions only if $\frac{d_1}{c_1}=\frac{d_2}{c_2}$; if that is the case, Solution set $S:(x,y,z)=\left(m,\ n,\frac{d_1}{c_1}\right)$, where $m,n\in\mathbb{Z}$ (lattice points on the plane $z=\frac{d_1}{c_1}$).

ii) *Group* 2: $\begin{bmatrix}0 & 0 & c_1\\ a_2 & 0 & 0\end{bmatrix}\begin{bmatrix}a_1 & 0 & 0\\ 0 & 0 & c_2\end{bmatrix}\begin{bmatrix}0 & 0 & c_1\\ a_2 & 0 & 0\end{bmatrix}$

$\begin{bmatrix}0 & 0 & c_1\\ a_2 & 0 & 0\end{bmatrix}$: The first equation described a plane parallel to or coincident with the *xy*-plane, while the second equation describes a plane parallel to or coincident with *yz*-plane. There are solutions only if $c_1\mid d_1$ and $a_2\mid d_2$ in which case, solution set $S:(x,y,z)=\left(\frac{d_2}{a_2},m,\frac{c_1}{d_1}\right)$; $m\in\mathbb{Z}$

iii) *Group* 3:

$\begin{bmatrix}0 & 0 & 0\\ 0 & b_2 & c_2\end{bmatrix}\begin{bmatrix}0 & 0 & 0\\ a_2 & 0 & c_2\end{bmatrix}\begin{bmatrix}0 & 0 & 0\\ a_2 & b_2 & 0\end{bmatrix}\begin{bmatrix}0 & b_1 & c_1\\ 0 & 0 & 0\end{bmatrix}\begin{bmatrix}a_1 & 0 & c_1\\ 0 & 0 & 0\end{bmatrix}\begin{bmatrix}a_1 & b_1 & 0\\ 0 & 0 & 0\end{bmatrix}$

$\begin{bmatrix} 0 & 0 & 0 \\ 0 & b_2 & c_2 \end{bmatrix}$: If $d_1 \neq 0$ there are no solutions. If $d_1 = 0$, the case reduces to part B(i).

iv) Exactly five of the six coefficients are zero

$$\begin{bmatrix} 0 & 0 & 0 \\ 0 & 0 & c_2 \end{bmatrix}\begin{bmatrix} 0 & 0 & 0 \\ 0 & b_2 & 0 \end{bmatrix}\begin{bmatrix} 0 & 0 & 0 \\ a_2 & 0 & 0 \end{bmatrix}\begin{bmatrix} 0 & 0 & c_1 \\ 0 & 0 & 0 \end{bmatrix}\begin{bmatrix} 0 & b_1 & 0 \\ 0 & 0 & 0 \end{bmatrix}\begin{bmatrix} a_1 & 0 & 0 \\ 0 & 0 & 0 \end{bmatrix}$$

$\begin{bmatrix} 0 & 0 & 0 \\ 0 & 0 & c_2 \end{bmatrix}$: If $d_1 \neq 0$ there are no solutions. If $d_1 = 0$, the case reduces to part B(ii).

5. **All coefficients are zero**

$\begin{bmatrix} 0 & 0 & 0 \\ 0 & 0 & 0 \end{bmatrix}$: If either $d_1$ or $d_2$ is non zero, there are no solutions. If $d_1 = d_2 = 0$, solution set $S:(x,y,z) = (m,n,k); m,n,k$ in other words each lattice point in space is a solution.